\documentclass[12pt,a4paper,english,intoc,bibliography=totoc,index=totoc,BCOR10mm,captions=tableheading,titlepage,fleqn]{scrbook}
\usepackage{lmodern}

\usepackage[T1]{fontenc}
\usepackage[latin9]{inputenc}
\usepackage{babel}
\usepackage{amsmath}
\usepackage{amssymb}
\usepackage[unicode=true,
 bookmarks=true,bookmarksnumbered=true,bookmarksopen=true,bookmarksopenlevel=2,
 breaklinks=false,pdfborder={0 0 0},backref=false,colorlinks=false]
 {hyperref}
\hypersetup{pdftitle={Abstract},
 pdfauthor={Uwe Stöhr},
 pdfpagelayout=OneColumn, pdfnewwindow=true, pdfstartview=XYZ, plainpages=false}
\usepackage{breakurl}

\makeatletter

\usepackage{calc}

\makeatother

\begin{document}
\addcontentsline{toc}{chapter}{Abstract} 

\noindent \begin{center}
\textbf{\Large Maximal independent sets in Borel graphs and large
cardinals}
\par\end{center}{\Large \par}

\noindent \begin{center}
{\large Haim Horowitz and Saharon Shelah}
\par\end{center}{\large \par}

\noindent \begin{center}
\textbf{Abstract}
\par\end{center}

\noindent \begin{center}
{\small We construct a Borel graph $G$ such that $ZF+DC+"$There
are no maximal independent sets in $G"$ is equiconsistent with $ZFC+"$There
exists an inaccessible cardinal$"$.}%
\footnote{{\small Date: June 14, 2016.}{\small \par}



Publication 1093 of the second author.%
}
\par\end{center}{\small \par}

\textbf{\large Introduction}{\large \par}

The main result of this note is motivated by our recent study of maximal
almost disjoint families and their relatives. Recall that $\mathcal F \subseteq [\omega]^{\omega}$
is a MAD family if $A \neq B \in \mathcal F \rightarrow |A \cap B|<\aleph_0$,
and $\mathcal F$ is maximal with respect to this property. Maximal
eventually different (MED) families are the analog of MAD families
where elements of $[\omega]^{\omega}$ are replaced by graphs of functions
from $\omega$ to $\omega$, namely, $f,g \in \omega^{\omega}$ are
eventually different if $f(n) \neq g(n)$ for large enough $n$, and
$\mathcal F \subseteq \omega^{\omega}$ is a MED family if the elements
of $\mathcal F$ are pairwise eventually different and $\mathcal F$
is maximal with respect to this property.

Questions on the (non-)existence and definability of such families
have attracted considerable interest for decades. The first results
were obtained by Mathias who proved the following theorem:

\textbf{Theorem {[}Ma{]}: }There are no analytic MAD families.

As for the possibility of the non-existence of MAD families, the following
result was recently proved by the authors (earlier such results were
proven by Mathias in {[}Ma{]} and by Toernquist in {[}To{]} using
Mahlo and inaccessible cardinals, respectively):

\textbf{Theorem {[}HwSh:1090{]}: }$ZF+DC+"$There are no MAD families$"$
is equiconsistent with $ZFC$.

Quite surprisingly, the situation for MED families turns out to be
different:

\textbf{Theorem {[}HwSh:1089{]}: }Assuming $ZF$, there exists a Borel
MED family.

A possible approach to explaining the above difference is via Borel
combinatorics. The study of Borel and analytic graphs was initiated
by Kechris, Solecki and Todorcevic in {[}KST{]}, and has been a source
of fruitful research ever since (see {[}KM{]} for a survey of recent
results). The above questions on MAD families are connected to Borel
combinatorics due to the following observation: There exist Borel
graphs $G_{MAD}$ and $H_{MED}$ such that there exists a MAD (MED)
family iff there exists a maximal independent set in $G_{MAD}$ $(H_{MED})$.
Therefore, we might try to explain the above difference of MAD and
MED families by pursuing the general problem of classifying Borel
graphs according to the consistency strength of $ZF+DC+"$There are
no maximal independent sets in $G"$.

The main goal of this note is to show that for some Borel graphs $G$,
$ZF+DC+"$There are no maximal independent sets in $G"$ has large
cardinal strength.

\textbf{\large The main result}{\large \par}

\textbf{Definition 1: }We shall define a Borel graph $G=(V,E)$ as
follows:

a. $V$ is the set of reals $r$ that code the following objects:

1. A linear order $I_r$ of the element of $\omega$ or some $n<\omega$.

2. A sequence $(s_{r,\alpha} : \alpha \in I_r)$ of pairwise distinct
reals.

3. A sequence of functions $(f_{r,a} : a \in I_r)$ such that each
$f_{r,a}$ is an injective function from $I_{r,<a}:= \{b \in I_r : b<_{I_r} a\}$
onto some initial segment of $\omega$.

b. Given $r_1 \neq r_2 \in V$ and $b\in I_{r_2}$, let $X_{r_1,r_2,b}$
be the set of pairs $(a_1,a_2) \in I_{r_1} \times I_{r_2,<b}$ such
that $s_{r_1,a_1}=s_{r_2,a_2}$.

c. Given $r_1 \neq r_2 \in V$, $\neg(r_1 E r_2)$ holds iff one of
the following holds:

1. There exists $b\in I_{r_2}$ such that $X_{r_1,r_2,b}$ is an isomorphism
from $I_{r_1}$ to $I_{r_2,<b}$ which also commutes with $f_{-,-}$.

2. There exists $b\in I_{r_1}$ such that $X_{r_2,r_1,b}$ is an isomorphism
from $I_{r_2}$ to $I_{r_1,<b}$ which also commutes with $f_{-,-}$.

\textbf{Definition 2: }Given $r_1 \neq r_2 \in V$, we say that $r_2$
extends $r_1$ and denote it by $r_1 <_G r_2$ when $\neg(r_1 E r_2)$
and clause (1) holds in definition 1(c).

\textbf{Claim 3 $(ZF+DC)$: }Let $X\subseteq V$ be an independent
set.

a. $X$ is linearly ordered by $<_G$.

b. If $X$ is countable then $X$ is not a maximal independent set.

\textbf{Proof: }a. Obvious.

b. By clause (a), there is a linear order $I$ such that $X=\{r_i : i\in I\}$
and $i<_I j$ iff $r_i <_G r_j$. For every $i<j \in I$, let $F_{i,j}$
be the isomorphism from $I_{r_i}$ to a proper initial segment of
$I_{r_j}$ witnessing $r_i <_G r_j$. Let $I_r$ be the direct limit
of the system $(I_{r_i},F_{j,k} : i,j,k \in I, j<k)$. For $a\in I_r$,
let $s_{r,a}$ be $s_{r_i,a'}$ where $a' \in I_{r_i}$ is some representative
of $a$, and define $f_{r,a}$ similarly. Let $r \in V$ be a real
coding $I_r$, $(s_{r,a} : a\in I_r)$ and $(f_{r,a} : a\in I_r)$,
then $\neg (rEr_i)$ for every $r_i \in X$. $\square$

\textbf{Theorem 4: }$ZF+DC+"$There is no maximal independent set
in $G"$ is equiconsistent with $ZFC+"$There exists an inaccessible
cardinal$"$.

Theorem 4 will follow from the following claims:

\textbf{Claim 5 $(ZF+DC)$: }If there exists $a\in \omega^{\omega}$
such that $\aleph_1=\aleph_1^{L[a]}$, then there exists a maximal
independent set in $G$.

\textbf{Claim 6: }There is no maximal independent set in $G$ in Levy's
model (aka Solovay's model).

Remark: While the set of vertices of $G$ is denoted by $V$, the
set-theoretic universe will be denoted by $\bold V$.

\textbf{Proof of claim 5: }Let $(s_{\alpha} : \alpha<\omega_1^{L[a]}) \in L[a]$
be a sequence of pairwise distinct reals, and let $\bar{f}^*=(f_{\alpha}^* : \alpha<\omega_1^{L[a]}) \in L[a]$
be a sequence of functions such that each $f_{\alpha}^*$ is an injective
function from $\alpha$ onto $|\alpha| \leq \omega$. For each $\alpha<\omega_1^{L[a]}$,
let $r_{\alpha} \in {(\omega^{\omega})}^{L[a]}$ be the $<_{L[a]}-$first
real that codes $(\alpha,(s_{\beta} : \beta<\alpha), \bar{f}^* \restriction \alpha)$.
The sequences $(s_{\alpha} : \alpha<\omega_1^{L[a]})$, $(f_{\alpha}^* : \alpha<\omega_1^{L[a]})$
and $(r_{\alpha} : \alpha<\omega_1^{L[a]})$ belong to $\bold V$,
and as $\omega_1=\omega_1^{L[a]}$, their length is $\omega_1$.

It's easy to see that $\{r_{\alpha} : \alpha<\omega_1^{L[a]}\}$ is
a well-defined set and is an independent subset of $V$, we shall
prove that it's a maximal independent set. Let $r\in V \setminus \{r_{\alpha} : \alpha<\omega_1^{L[a]}\}$
and suppose towards contradiction that $\neg(rEr_{\alpha})$ for every
$\alpha<\omega_1^{L[a]}$. There are two possible cases:

\textbf{Case I: }$r_{\alpha}<_G r$ for every $\alpha<\omega_1^{L[a]}$.
In this case, $I_r$ is a linear order, and each $\alpha<\omega_1^{L[a]}$
embeds into $I_r$ as an initial segment. Therefore, $\omega_1=\omega_1^{L[a]}$
embeds into $I_r$ as an initial segment, a contradiction.

\textbf{Case II: }$r<_G r_{\alpha}$ for some $\alpha<\omega_1^{L[a]}$.
Let $\alpha$ be the minimal ordinal with this property, then $\alpha$
necessarily has the form $\beta+1$. If $r=r_{\beta}$, then we get
a contradiction to the choice of $r$. If $r \neq r_{\beta}$, then
it's easy to see that $rEr_{\beta}$, contradicting our assumption.
$\square$

\textbf{Proof of claim 6: }Let $\kappa$ be an inaccessible cardinal
and let $\mathbb P=Coll(\aleph_0,<\kappa)$, we shall prove that $\Vdash_{\mathbb P} "$There
is no maximal independent set in $G$ from $HOD(\mathbb R)"$. Suppose
towards contradiction that $p\in \mathbb P$ forces that $\underset{\sim}{X}$
is such a set. Let $\mathbb Q$ be a forcing notion such that $\mathbb Q \lessdot \mathbb P$,
$|\mathbb Q|<\kappa$, $p\in \mathbb Q$ and $\underset{\sim}{X}$
is definable using a parameter from $\mathbb{R}^{{\bold V}^{\mathbb Q}}$.
By the properties of the Levy collapse, we may assume wlog that $\mathbb Q=\{0\}$
and $p=0$. If $\Vdash_{\mathbb P} "\underset{\sim}{X} \subseteq (\omega^{\omega})^{\bold V}"$,
then $\Vdash_{\mathbb P} "|\underset{\sim}{X}|=\aleph_0"$, and by
claim 3, $\underset{\sim}{X}$ is not a maximal independent set in
${\bold V}^{\mathbb P}$, a contradiction. Therefore, there exist
$p_1 \in \mathbb P$ and $\underset{\sim}{r_1}$ such that $p_1 \Vdash_{\mathbb P} "\underset{\sim}{r_1} \in \underset{\sim}{X} \wedge \underset{\sim}{r_1} \notin \bold V"$.
Let $\mathbb{Q}_1 \lessdot \mathbb P$ be a forcing of cardinality
$<\kappa$ such that $p_1 \in \mathbb{Q}_1$ and $\underset{\sim}{r_1}$
is a $\mathbb{Q}_1$-name. For $l=2,3$ let $(\mathbb{Q}_l,p_l,\underset{\sim}{r_l})$
be isomorphic copies of $(\mathbb{Q}_1,p_1,\underset{\sim}{r_1})$
such that $\underset{n=1,2,3}{\Pi}\mathbb{Q}_n \lessdot \mathbb P$
(identifying $\mathbb{Q}_1$ with its canonical image in the product).
Choose $(p_1,p_2) \leq (q_1,q_2)$ such that $(q_1,q_2) \Vdash_{\mathbb{Q}_1 \times \mathbb{Q}_2} "\underset{\sim}{r_1} \neq \underset{\sim}{r_2}"$.
As $(q_1,q_2) \Vdash_{\mathbb{Q}_1 \times \mathbb{Q}_2} "\underset{\sim}{r_1},\underset{\sim}{r_2} \in \underset{\sim}{X}"$,
then wlog $(q_1,q_2)$ forces that $\underset{\sim}{r_1}<_G \underset{\sim}{r_2}$
as witnessed by an isomorphism from $I_{\underset{\sim}{r_1}}$ to
$I_{\underset{\sim}{r_2},<s}$ for some $s\in I_{\underset{\sim}{r_2}}$.
Let $q_3 \in \mathbb{Q}_3$ be the conjugate of $q_1$, then $(q_2,q_3)$
forces (in $\mathbb{Q}_2 \times \mathbb{Q}_3$) that $\underset{\sim}{r_2},\underset{\sim}{r_3} \in \underset{\sim}{X}$
and $\underset{\sim}{r_3}<_G \underset{\sim}{r_2}$ as witnessed by
an isomorphism from $I_{\underset{\sim}{r_3}}$ to $I_{\underset{\sim}{r_2},<s}$.
Now pick $(q_1,q_2,q_3)\leq (q_1',q_2',q_3')$ that forces in addition
that $\underset{\sim}{r_1} \neq \underset{\sim}{r_3}$, then necessarily
it forces that $\underset{\sim}{r_1}E \underset{\sim}{r_3}$, a contradiction.
$\square$

\textbf{\large Open problems}{\large \par}

\textbf{Notation: }Given a Borel graph $G$, let $\psi(G)$ be the
statement $"$There are no maximal independent sets in $G"$.

\textbf{Problem 1: }Classify the Borel graphs according to the consistency
strength of $ZF+DC+\psi(G)$.

As the above problem seems to be quite difficult at the moment, it
might be reasonable to consider the following subproblems first:

\textbf{Problem 2: }What are the possibilites (in terms of large cardinal
strength) for the consistency strength of $ZF+DC+\psi(G)$?

\textbf{Problem 3: }Find combinatorial/descriptive set theoretic/model
theoretic properties $\phi_1$ and $\phi_2$ such that:

a. $\phi_1(G_{MAD})$.

b. $\phi_2(H_{MED})$.

c. $\phi_1(G) \rightarrow ZF+DC+\psi(G)$ is equiconsistent with $ZFC$.

d. $\phi_2(G) \rightarrow ZF+DC \vdash \neg \psi(G)$.

e. $\phi_1$ and $\phi_2$ are satisfied by a large collection of
Borel graphs.

A solution to problem (3) would explain the difference between MAD
and MED families that was discussed in the introduction.

\textbf{\large References}{\large \par}

{[}HwSh:1089{]} Haim Horowitz and Saharon Shelah, A Borel maximal
eventually different family, arXiv:1605.07123.

{[}HwSh:1090{]} Haim Horowitz and Saharon Shelah, Can you take Toernquist's
inaccessible away?, arXiv:1605.02419.

{[}KM{]} A. S. Kechris and A. Marks, Descriptive graph combinatorics,
Manuscript, 2015.

{[}KST{]} A. S. Kechris, S. Solecki and S. Todorcevic, Borel chromatic
numbers, Adv. Math. \textbf{141 }(1999), 1-44.

{[}Ma{]} A. R. D Mathias, Happy families, Ann. Math. Logic \textbf{12
}(1977), no. 1, 59-111. MR 0491197. 

{[}To{]} Asger Toernquist, Definability and almost disjoint families,
arXiv:1503.07577.

$\\$

(Haim Horowitz) Einstein Institute of Mathematics

Edmond J. Safra campus,

The Hebrew University of Jerusalem.

Givat Ram, Jerusalem, 91904, Israel.

E-mail address: haim.horowitz@mail.huji.ac.il

$\\$

(Saharon Shelah) Einstein Institute of Mathematics

Edmond J. Safra campus,

The Hebrew University of Jerusalem.

Givat Ram, Jerusalem, 91904, Israel.

Department of Mathematics

Hill Center - Busch Campus,

Rutgers, The State University of New Jersey.

110 Frelinghuysen road, Piscataway, NJ 08854-8019 USA

E-mail address: shelah@math.huji.ac.il
\end{document}